\documentclass[11pt]{article}
\usepackage{mathrsfs,amsmath,amssymb,amsfonts}
\numberwithin{equation}{section}
\begin{document}
\title{Pointwise weighted approximation of functions with inner
singularities by combinations of Bernstein operators}
\author{Wen-ming Lu$^{1}$, Lin Zhang$^{2}$\thanks{Corresponding author. Address: Department of Mathematics, Zhejiang University, Hangzhou 310027, PR China. \textsl{E-mail address}:
\texttt{linyz@zju.edu.cn}(L.-Zhang);\texttt{lu\_wenming@163.com}(W.-Lu).},
and Meng-yi Chai$^{3}$}
\date{\it $^{1,3}$School of Science, Hangzhou Dianzi University,
Hangzhou, 310018 P.R. China\\$^2$Department of Mathematics, Zhejiang
University, Hangzhou, 310027 P.R. China} \maketitle
\mbox{}\hrule\mbox{}\\[0.5cm]
\textbf{Abstract}\\[-0.2cm]^^L
We introduce another new type of combinations of Bernstein operators in this paper, which can be used to
approximate the functions with inner singularities.
The direct and inverse results of the weighted approximation of this new type of combinations are obtained.\\~\\
%---------------------------------------------------------------------------------------------------------
\mbox{}\hrule\mbox{}
%------------------------------------------------------------------------------------------
\section{Introduction}
The set of all continuous functions, defined on the interval $I$, is
denoted by $C(I)$. For any $f\in C([0,1])$, the corresponding
\emph{Bernstein operators} are defined as follows:
$$B_n(f,x):=\sum_{k=0}^nf(\frac{k}{n})p_{n,k}(x),$$
where
$$p_{n,k}(x):={n \choose k}x^k(1-x)^{n-k}, \ k=0,1,2,\ldots,n, \ x\in[0,1].$$
Approximation properties of Bernstein operators have been studied
very well (see \cite{Berens}, \cite{Della},
\cite{Totik}-\cite{Lorentz}, \cite{Yu}-\cite{X. Zhou}, for example).
In order to approximate the functions with singularities, Della
Vecchia et al. \cite{Della} and Yu-Zhao \cite{Yu} introduced some
kinds of modified Bernstein operators. Throughout the paper,
$C$ denotes a positive constant independent of $n$ and $x$,
which may be different in different cases.\\~\\
%--------------------------------------------------------------------------------------------
Let $\bar{w}(x)=|x-\xi|^\alpha,\ 0<\xi<1,\ \alpha>0$ and
$C_{\bar{w}}:= \{{f \in C([0,1] \setminus {\xi})
:\lim\limits_{x\longrightarrow\xi}(\bar{w}f)(x)=0 }\}.$ The norm in
$C_{\bar{w}}$ is defined as
$\|f\|_{C_{\bar{w}}}:=\|{\bar{w}}f\|=\sup\limits_{0\leqslant
x\leqslant 1}|({\bar{w}f})(x)|$. Define
$$W_{\bar{w},\lambda}^{r}:= \{f \in
C_{\bar{w}}:f^{(r-1)} \in A.C.((0,1)),\
\|{\bar{w}}\varphi^{r\lambda}f^{(r)}\|<\infty\}.$$ For $f \in
C_{\bar{w}}$, define the \emph{weighted modulus of smoothness}
by\\
$$\omega_{\varphi^\lambda}^{r}(f,t)_{\bar{w}}:=\sup_{0<h\leqslant
t}\{\|{\bar{w}}\triangle_{h\varphi^\lambda}^{r}f\|_{[16h^2,1-16h^2]}+\|{\bar{w}}{\overrightarrow{\triangle}_{h}^{r}f\|_{[0,16h^2]}+\|{\bar{w}}{\overleftarrow{\triangle}_{h}^{r}}f\|_{[1-16h^2,1]}}\},$$
where
\begin{eqnarray*}
\Delta_{h\varphi^\lambda}^{r}f(x)&=&\sum_{k=0}^{r}(-1)^{k}{r \choose
k}f(x+(\frac r2-k)h\varphi^\lambda(x)),\\
\overrightarrow{\Delta}_{h}^{r}f(x)&=&\sum_{k=0}^{r}(-1)^{k}{r \choose k}f(x+(r-k)h),
\end{eqnarray*}
and $\varphi(x)=\sqrt{x(1-x)}$.  On the other hand,
since the Bernstein polynomials cannot be used for the
investigation of higher orders of smoothness, Butzer \cite{Butzer}
introduced the combinations of Bernstein polynomials which
have higher orders of approximation. Ditzian and Totik \cite{Totik}
extended this method of combinations and defined the following
combinations of Bernstein operators:
\begin{eqnarray*}
B_{n,r}(f,x):=\sum_{i=0}^{r-1}C_{i}(n)B_{n_i}(f,x).
\end{eqnarray*}
With the conditions\\
(a) $n=n_0<n_1< \cdots <n_{r-1}\leqslant Cn,$\\
(b) $\sum_{i=0}^{r-1}|C_{i}(n)|\leqslant C,$\\
(c) $\sum_{i=0}^{r-1}C_{i}(n)=1,$\\
(d) $\sum_{i=0}^{r-1}C_{i}(n)n_{i}^{-k}=0$, for $k=1,\cdots,r-1$.
%------------------------------------------------------------------------------------------
\section{The main results}
For any positive integer $r$, we consider the determinant
\begin{eqnarray*}
A_{r}:=
\begin {matrix}
\begin{vmatrix}
1 & 1 & 1 & \cdots & 1 \\
2r+1 & 2r+2 & 2r+3 & \cdots & 4r+1 \\
(2r)(2r+1) & (2r+1)(2r+2) & (2r+2)(2r+3) & \cdots & (4r)(4r+1) \\
\cdots & \cdots & \cdots & \ddots & \cdots \\
2\cdots(2r+1) & 3\cdots(2r+2) & 4\cdots(2r+3) & \cdots &
(2r+2)\cdots(4r+1) \end{vmatrix} &
\end{matrix}.
\end{eqnarray*}
We obtain $A_{r}=\prod_{j=2}^{2r}j!$. Thus, there is a
unique solution for the system of nonhomogeneous linear equations:
\begin{eqnarray}
\left\{
   \begin{array}{ccccccccc}
   a_1 & + & a_2 & + & \cdots & + & a_{2r+1} & = &1, \\
   (2r+1)a_1 & + & (2r+2)a_2 & + & \cdots & + & (4r+1)a_{2r+1} & = &0, \\
   (2r+1)(2r)a_1 & + & (2r+1)(2r+2)a_2 & + & \cdots & + & (4r)(4r+1)a_{2r+1} & = &0, \\
   &&&\vdots&&&  \\
   (2r+1)!a_1 & + & 3 \cdots (2r+2)a_2 & + & \cdots & + & (2r+2) \cdots (4r+1)a_{2r+1} & =
   &0.
   \end{array} \right.\label{s1}
\end{eqnarray}
Let
\begin{eqnarray*}
\psi(x)=\left\{
\begin{array}{lrr}
a_1x^{2r+1}+a_2x^{2r+2}+\cdots+a_{2r+1}x^{4r+1}, &&0<x<1, \\
0,   &&x \leqslant0,  \\
1,  &&x=1.
             \end{array}
\right.
\end{eqnarray*}
with the coefficients $a_1,$ $a_2,$ $\cdots,$ $a_{2r+1}$ satisfying
(\ref{s1}). From (\ref{s1}), we see that $\psi(x)\in
C^{(2r)}(-\infty,+\infty)$, $0\leqslant\psi(x)\leqslant1$ for $0
\leqslant x\leqslant1$. Moreover, it holds that $\psi(1)=1$,\
$\psi^{(i)}(0)=0,\ i=0,1,\cdots,2r$ and $\psi^{(i)}(1)=0,\
i=1,2,\cdots,2r$. \\
Let $$H(f,x):=\sum_{i=1}^{r+1}f(x_{i})l_{i}(x),$$ and
$$l_{i}(x):=\frac{\prod_{j=1,j\neq i}^{r+1}(x-x_{j})}{\prod_{j=1,j\neq i}^{r+1}(x_{i}-x_{j})},\ x_{i}=\frac{[n\xi-({(r-1)/2}+i)]}{n},\ i=1,2, \cdots r+1.$$
Further, let
$$x_{1}^{'}=\frac{[n\xi-2\sqrt{n}]}{n},\ x_{2}^{'}=\frac{[n\xi-\sqrt{n}]}{n},\ x_{3}^{'}=\frac{[n\xi+\sqrt{n}]}{n},\ x_{4}^{'}=\frac{[n\xi+2\sqrt{n}]}{n},$$
and
$${\bar{\psi}}_{1}(x)=\psi(\frac{x-x_{1}^{'}}{x_{2}^{'}-x_{1}^{'}}),\ {\bar{\psi}}_{2}(x)=\psi(\frac{x-x_{3}^{'}}{x_{4}^{'}-x_{3}^{'}}).$$
Set
$${\bar{F}}_{n}(f,x):={\bar{F}}_{n}(x)=f(x)(1-{\bar{\psi}}_{1}(x)+{\bar{\psi}}_{2}(x))+{\bar{\psi}}_{1}(x)(1-{\bar{\psi}}_{2}(x))H(x).$$
We have
\begin{eqnarray*}
{\bar{F}}_{n}(f,x)=\left\{\begin{array}{lr}
f(x),          &       x\in [0,x_{r-5/2}]\cup [x_{r+3/2},1],   \\
f(x)(1-{\bar{\psi}}_{1}(x))+{\bar{\psi}}_{1}(x)
H(x),      &
x\in [x_{r-5/2},x_{r-3/2}],  \\
H(x),          &       x\in [x_{r-3/2},x_{r+1/2}],  \\
H(x)(1-{\bar{\psi}}_{2}(x))+{\bar{\psi}}_{2}(x)f(x), & x\in
[x_{r+1/2},x_{r+3/2}].
            \end{array}
\right.
\end{eqnarray*}
Obviously, ${\bar{F}}_{n}(f,x)$ is linear, reproduces
polynomials of degree $r$, and ${\bar{F}}_{n}(f,x)\in
C^{(2r)}([0,1])$, provided that $f \in C^{(2r)}([0,1]).$
%----------------------------------------------------------------------------------
Now, we can define our new combinations of Bernstein operators as
follows:
\begin{eqnarray}
{\bar{B}}_{n,r}(f,x):=B_{n,r}({\bar{F}_{n}},x)=\sum_{i=0}^{r-1}C_{i}(n)B_{n_i}({\bar{F}_{n}},x),\label{s2}
\end{eqnarray}
where $C_{i}(n)$ satisfy the conditions (a)-(d). Our main result is
the following: \\~\\
\textbf{Theorem.} \textit{For $f\in C_{\bar{w}},\ 0 \leqslant
\lambda \leqslant 1,\   0<\xi<1,\ \alpha>0,\ 0 < \alpha_0 <r,$ we
have
\begin{eqnarray*}
{\bar{w}(x)}|f(x)-{\bar{B}_{n,r-1}(f,x)}|=O((n^{-{\frac
12}}\varphi^{-\lambda}(x)\delta_n(x))^{\alpha_0})
\Longleftrightarrow
\omega_{\varphi^\lambda}^r(f,t)_{\bar{w}}=O(t^{\alpha_0}).
\end{eqnarray*}}
%------------------------------------------------------------------------------------------
\section{Lemmas}
\textbf{Lemma 1.}(\cite{Della}) \textit{If $\gamma \in R,$ then
\begin{eqnarray}
\sum_{k=0}^np_{n,k}(x)|k-nx|^\gamma \leqslant Cn^{\frac
\gamma2}\varphi^\gamma(x).\label{s3}
\end{eqnarray}}
\textbf{Lemma 2.}(\cite{Lu}) \textit{Let
$A_n(x):={\bar{w}(x)}\sum\limits_{|k-n\xi|\leqslant
\sqrt{n}}p_{n,k}(x)$. Then $A_n(x)\leqslant Cn^{-\alpha/2}$ for
$0<\xi <1$ and $\alpha>0$.} \\~\\
%----------------------------------------------------------------------
\textbf{Lemma 3.} \textit{For any $\alpha >0,$ $0\leqslant \lambda
\leqslant1,\ f\in C_{\bar{w}},$ we have
\begin{eqnarray}
\|{\bar{w}}{\bar{B}}_{n,r-1}^{(r)}(f)\|\leqslant
Cn^{r}\|{\bar{w}}f\|.\label{s4}
\end{eqnarray}}
\textbf{Proof.} We first prove $x\in [0,{\frac 1n})$ (The same as
$x\in (1-{\frac 1n},1]$), now
\begin{eqnarray*}
|{\bar{w}}(x){\bar{B}}_{n,r-1}^{(r)}(f,x)|&\leqslant&
{\bar{w}}(x)\sum_{i=0}^{r-2}{\frac
{n_{i}!}{({n_{i}-r})!}}\sum_{k=0}^{n_i-r}|C_{i}(n)||\overrightarrow{\Delta}_{\frac
1{n_i}}^{r}{\bar{F}}_{n}{(\frac k{n_i})}|p_{n_i-r,k}(x)\nonumber\\
&\leqslant&
C{\bar{w}}(x)\sum_{i=0}^{r-2}n_{i}^{r}\sum_{k=0}^{n_i-r}|\overrightarrow{\Delta}_{\frac
1{n_i}}^{r}{\bar{F}}_{n}{(\frac k{n_i})}|p_{n_i-r,k}(x)\nonumber\\
&\leqslant&
C{\bar{w}}(x)\sum_{i=0}^{r-2}n_{i}^{r}\sum_{k=0}^{n_i-r}\sum_{j=0}^{r}|{\bar{F}}_{n}({\frac
{k+r-j}{n_i}})|p_{n_i-r,k}(x)\nonumber\\
&\leqslant&
C{\bar{w}}(x)\sum_{i=0}^{r-2}n_{i}^{r}\sum_{j=0}^{r}|{\bar{F}}_{n}({\frac
{r-j}{n_i}})|p_{n_i-r,0}(x)\nonumber\\
&&+
C{\bar{w}}(x)\sum_{i=0}^{r-2}n_{i}^{r}\sum_{j=0}^{r}|{\bar{F}}_{n}({\frac
{n_{i}-j}{n_i}})|p_{n_i-r,n_i-r}(x)\nonumber\\
&&+
C{\bar{w}}(x)\sum_{i=0}^{r-2}n_{i}^{r}\sum_{k=1}^{n_i-r-1}\sum_{j=0}^{r}|{\bar{F}}_{n}({\frac
{k+r-j}{n_i}})|p_{n_i-r,k}(x)\nonumber\\
&:=&H_1 +H_2 + H_3.
\end{eqnarray*}
We have
\begin{eqnarray*}
H_1&\leqslant&
C{\bar{w}}(x)\sum_{i=0}^{r-2}n_{i}^{r}\sum_{j=0}^{r}|{\bar{F}}_{n}({\frac
{r-j}{n_i}})|p_{n_i-r,0}(x)\nonumber\\
&\leqslant&
Cn^{r}\|{\bar{w}}f\|\sum_{i=0}^{r-2}\sum_{j=0}^{r}(\frac{n_{i}|x-\xi|}{r-j-n_{i}\xi})^{\alpha}(1-x)^{{n_{i}}-r}\nonumber\\
&\leqslant&
Cn^{r}\|{\bar{w}}f\|\sum_{i=0}^{r-2}(n_{i}|x-\xi|)^{\alpha}(1-x)^{{n_{i}}-r}\nonumber\\
&\leqslant& Cn^{r}\|{\bar{w}}f\|.
\end{eqnarray*}
Similarly, we can get $H_2\leqslant Cn^{r}\|{\bar{w}}f\|,$
and $H_3\leqslant
Cn^{r}\|{\bar{w}}f\|$. \\~\\
%-----------------------------------------------------------------------------------------------------
When $x\in [{\frac 1n},1-{\frac 1n}],$ according to
\cite{Totik}, we have
\begin{eqnarray*}
&&|{\bar{w}}(x){\bar{B}}_{n,r-1}^{(r)}(f,x)|\nonumber\\
&=&|{\bar{w}}(x)B_{n,r-1}^{(r)}({\bar{F}_{n}},x)|\nonumber\\
&=&{\bar{w}}(x)(\varphi^{2}(x))^{-r}\sum_{i=0}^{r-2}\sum_{j=0}^{r}Q_{j}(x,n_i)|C_{i}(n)|n_{i}^{j}\sum_{k/n_i\in
A}|(x-{\frac kn_{i}})^{j}||{\bar{F}}_{n}({\frac kn_{i}})|p_{n_i,k}(x)\nonumber\\
&&+{\bar{w}}(x)(\varphi^{2}(x))^{-r}\sum_{i=0}^{r-2}\sum_{j=0}^{r}Q_{j}(x,n_i)|C_{i}(n)|n_{i}^{j}\sum_{x_2^{\prime}
\leqslant k/n_i\leqslant x_3^\prime}|(x-{\frac kn_{i}})^{j}||H({\frac
kn_{i}})|p_{n_i,k}(x)\nonumber\\
&:=&\sigma_1+ \sigma_2.
\end{eqnarray*}
Where $A:=[0,x_2^{\prime}]\cup [x_3^{\prime},1]$, $H$ is a
linear function. If ${\frac kn_{i}}\in A,$ when ${\frac
{\bar{w}(x)}{\bar{w}(\frac {k}{n_{i}})}}\leqslant C(1+n_i^{-{\frac
{\alpha}{2}}}|k-n_ix|^\alpha),$ we have $|k-n_{i}\xi|\geqslant
{\frac {\sqrt{n_{i}}}{2}}$, also
$Q_{j}(x,n_i)=(n_ix(1-x))^{[{\frac {r-j}{2}}]},$ and
$(\varphi^{2}(x))^{-r}Q_{j}(x,n_i)n_{i}^{j}\leqslant
C(n_i/\varphi^{2}(x))^{\frac {r+j}{2}}.$\\
By (\ref{s3}), then
\begin{eqnarray*}
\sigma_1&\leqslant&
C{\bar{w}}(x)\sum_{i=0}^{r-2}\sum_{j=0}^{r}(\frac
{n_{i}}{\varphi^{2}(x)})^{\frac {r+j}{2}}\sum_{k=0}^{n_{i}}|(x-{\frac
kn_{i}})^{j}||{\bar{F}}_{n}({\frac kn_{i}})|p_{n_i,k}(x)\nonumber\\
&\leqslant&
C\|{\bar{w}}f\|\sum_{i=0}^{r-2}\sum_{j=0}^{r}(\frac
{n_{i}}{\varphi^{2}(x)})^{\frac {r+j}{2}}\sum_{k=0}^{n_{i}}[1+n_i^{-{\frac
{\alpha}{2}}}|k-n_ix|^\alpha]|x-{\frac kn_{i}|^{j}}p_{n_i,k}(x)\nonumber\\
&:=&I_1 + I_2.
\end{eqnarray*}
By a simple calculation, we have $I_1\leqslant
Cn^{r}\|{\bar{w}}f\|$. By ({\ref{s3}), then
$$I_2\leqslant C\|{\bar{w}}f\|\sum_{i=0}^{r-2}\sum_{j=0}^{r}n_i^{-({{\frac {\alpha}{2}}}+j)}(\frac {n_{i}}
{\varphi^{2}(x)})^{\frac {j}{2}}\sum_{k=0}^{n_{i}}|k-n_ix|^{\alpha+j}p_{n_i,k}(x)\leqslant
Cn^{r}\|{\bar{w}}f\|.$$\\
We note that $|H({\frac {k}{n_i}})|\leqslant
max(|H(x_1^\prime)|,|H(x_4^\prime)|):=H(a)$.\\~\\
%----------------------------------------------------------------------------------------------------------------
If $x\in [x_1^\prime,x_4^\prime],$ we have
${\bar{w}(x)}\leqslant {\bar{w}(a)}.$ So, if $x\in
[x_1^\prime,x_4^\prime],$ then
$$\sigma_2\leqslant Cn^r{\bar{w}(a)}H(a)\leqslant Cn^{r}\|{\bar{w}}f\|.$$
%----------------------------------------------------------------------------------------------------------------
If $x\notin [x_1^\prime,x_4^\prime],$ then
${\bar{w}(a)}>n_i^{-{\frac {\alpha}{2}}},$ by lemma 2, we have
$$\sigma_2\leqslant C{\bar{w}(a)}H(a){\bar{w}}(x)\sum_{i=0}^{r-2}n_{i}^{r+{\frac {\alpha}{2}}}
\sum_{x_2^{\prime} \leqslant k/n_i\leqslant
x_3^\prime}p_{n_i,k}(x)\leqslant Cn^{r}\|{\bar{w}}f\|.$$\\
It follows from combining the above inequalities that the lemma is
proved. $\Box$\\~\\
\textbf{Lemma 4.}(\cite{Lu}) \textit{For any $\alpha >0,$
$0\leqslant \lambda \leqslant1,\ f\in C_{\bar{w}},$ we have
\begin{eqnarray}
\|{\bar{w}}{\bar{B}}_{n,r-1}(f)\|\leqslant
C\|{\bar{w}}f\|.\label{s5}
\end{eqnarray}}
%-------------------------------------------------------------------------------------------------------
\textbf{Lemma 5.} (\cite{Wang}) \textit{If \
$\varphi(x)=\sqrt{x(1-x)},\ 0 \leqslant \lambda \leqslant 1,\ 0
\leqslant \beta \leqslant 1,\ \alpha>0,$ then
\begin{eqnarray}
\int_{-{\frac {h\varphi^\lambda(x)}{2}}}^{\frac
{h\varphi^\lambda(x)}{2}} \cdots \int_{-{\frac
{h\varphi^\lambda(x)}{2}}}^{\frac
{h\varphi^\lambda(x)}{2}}\varphi^{-r\beta}(x+\sum_{k=1}^ru_k)du_1
\cdots du_r \leqslant Ch^r\varphi^{r(\lambda-\beta)}(x).\label{s6}
\end{eqnarray}}
\textbf{Lemma 6.} \textit{For any $r \in N,\ f\in
W_{\bar{w},\lambda}^{r},\ 0 \leqslant \lambda \leqslant 1,\
\alpha>0,$ we have
\begin{eqnarray}
\|{\bar{w}}\varphi^{r\lambda}{\bar{F}_n}^{(r)}\| \leqslant
C\|{\bar{w}}\varphi^{r\lambda}f^{(r)}\|.\label{s7}
\end{eqnarray}}
\textbf{Proof.} We first prove $x\in [x_{r-5/2},x_{r-3/2}]$ (The
same as the others), we have
\begin{eqnarray*}
|{\bar{w}}(x)\varphi^{r\lambda}(x){\bar{F}_n}^{(r)}(x)| &\leqslant&
|{\bar{w}}(x)\varphi^{r\lambda}(x)f^{(r)}(x)| +
|{\bar{w}}(x)\varphi^{r\lambda}(x)(f(x)-{\bar{F}_n}(x))^{(r)}|\\
&:=& I_1 + I_2.
\end{eqnarray*}
Obviously
\begin{eqnarray*}
I_1 \leqslant C\|{\bar{w}}\varphi^{r\lambda}f^{(r)}\|.
\end{eqnarray*}
For $I_2,$ we have
\begin{eqnarray*}
I_2 =
{\bar{w}}(x)\varphi^{r\lambda}(x)|(f(x)-{\bar{F}_n}(x))^{(r)}|={\bar{w}}(x)\varphi^{r\lambda}(x)\sum_{i=0}^rn^{\frac
i2}|(f(x)-{\bar{F}_n}(x))^{(r-i)}|.
\end{eqnarray*}
By \cite{Totik}, we have
\begin{eqnarray*}
|(f(x)-{\bar{F}_n}(x))^{(r-i)}|_{[x_{r-5/2},x_{r-3/2}]} \leqslant
C(n^{(r-i)/2}\|f-H\|_{[x_{r-5/2},x_{r-3/2}]} + n^{-
i/2}\|f^{(r)}\|_{[x_{r-5/2},x_{r-3/2}]}).
\end{eqnarray*}
So
\begin{eqnarray*}
I_2 &\leqslant& Cn^{\frac
r2}{\bar{w}}(x)\varphi^{r\lambda}(x)\|f-H\|_{[x_{r-5/2},x_{r-3/2}]}
+
C{\bar{w}}(x)\varphi^{r\lambda}(x)\|f^{(r)}\|_{[x_{r-5/2},x_{r-3/2}]}\\
&:=& T_1 + T_2.
\end{eqnarray*}
By Taylor expansion, we have
\begin{eqnarray}
f({x_i})=\sum_{u=0}^{r-1}\frac{(x_i-x)^u}{u!}f^{(u)}(x)+{\frac
1{(r-1)!}}\int_{x}^{x_{i}}(x_i-s)^{r-1}f^{(r)}(s)ds,\label{s8}
\end{eqnarray}
It follows from (\ref{s8}) and the identity
\begin{eqnarray*}
\sum\limits_{i=1}^{r}x_{i}^{v}l_{i}(x)=Cx^v,\ v=0,1,\cdots,r.
\end{eqnarray*}
we have
\begin{eqnarray*}
H(f,x)&=&\sum_{i=1}^{r}\sum_{u=0}^{r}\frac{(x_i-x)^u}{u!}f^{(u)}(x)l_{i}(x)+{\frac
1{(r-1)!}}\sum_{i=1}^{r}l_{i}(x)\int_{x}^{x_{i}}(x_i-s)^{r-1}f^{(r)}(s)ds\nonumber\\
&=&f(x)+C\sum_{u=1}^{r}f^{(u)}(x)(\sum_{v=0}^{u}C_{u}^{v}(-x)^{u-v}\sum_{i=1}^{r}x_{i}^{v}l_{i}(x))\nonumber\\
&&+{\frac
1{(r-1)!}}\sum_{i=1}^{r}l_{i}(x)\int_{x}^{x_{i}}(x_i-s)^{r-1}f^{(r)}(s)ds,
\end{eqnarray*}
which implies that
\begin{eqnarray*}
{\bar{w}(x)}\varphi^{r\lambda}(x)(f(x)-H(f,x))={\frac
1{(r-1)!}}{\bar{w}(x)}\varphi^{r\lambda}(x)\sum_{i=1}^{r}l_{i}(x)\int_{x}^{x_{i}}(x_i-s)^{r-1}f^{(r)}(s)ds,
\end{eqnarray*}
since $|l_{i}(x)|\leqslant C$ for $x\in [x_{r-5/2},x_{r-3/2}],\
i=1,2,\cdots,r$. It follows from
$\frac{|x_i-s|^{r-1}}{{\bar{w}}(s)}\leqslant
\frac{|x_i-x|^{r-1}}{{\bar{w}}(x)},$ $s$ between $x_i$ and $x$, then
\begin{eqnarray*}
{\bar{w}(x)}\varphi^{r\lambda}(x)|f(x)-H(f,x)|&\leqslant&C\bar{w}(x)\varphi^{r\lambda}(x)\sum_{i=1}^{r}\int_{x}^{x_{i}}(x_i-s)^{r-1}|f^{(r)}(s)|ds\nonumber\\
&\leqslant&C\varphi^{r\lambda}(x)\|{\bar{w}}\varphi^{r\lambda}f^{(r)}\|\sum_{i=1}^{r}(x_i-x)^{r-1}\int_{x}^{x_{i}}\varphi^{-r\lambda}(s)ds\nonumber\\
&\leqslant&{\frac
{C}{n^{r/2}}}\|{\bar{w}}\varphi^{r\lambda}f^{(r)}\|.
\end{eqnarray*}
So
\begin{eqnarray*}
I_2 \leqslant C\|{\bar{w}}\varphi^{r\lambda}f^{(r)}\|.
\end{eqnarray*}
Then, the lemma is proved. $\Box$\\~\\
\textbf{Lemma 7.} \textit{For any $g\in W_{\bar{w},\lambda}^{r},\ 0
\leqslant \lambda \leqslant 1,\ \alpha>0,$ we have
\begin{eqnarray}
{\bar{w}(x)}|g(x)-H(g,x)| \leqslant C(\frac
{\delta_n(x)}{\sqrt{n}\varphi^\lambda(x)})^r\|{\bar{w}}\varphi^{r\lambda}g^{(r)}\|.\label{s9}
\end{eqnarray}}
\textbf{Proof.} By Taylor expansion, we have
\begin{eqnarray*}
f({x_i})=\sum_{u=0}^{r-1}\frac{(x_i-x)^u}{u!}f^{(u)}(x)+{\frac
1{(r-1)!}}\int_{x}^{x_{i}}(x_i-s)^{r-1}f^{(r)}(s)ds,
\end{eqnarray*}
It follows from the above equality and the identity
\begin{eqnarray*}
\sum\limits_{i=1}^{r}x_{i}^{v}l_{i}(x)=Cx^v,\ v=0,1,\cdots,r.
\end{eqnarray*}
we have
\begin{eqnarray*}
H(f,x)&=&\sum_{i=1}^{r}\sum_{u=0}^{r}\frac{(x_i-x)^u}{u!}f^{(u)}(x)l_{i}(x)+{\frac
1{(r-1)!}}\sum_{i=1}^{r}l_{i}(x)\int_{x}^{x_{i}}(x_i-s)^{r-1}f^{(r)}(s)ds\nonumber\\
&=&f(x)+C\sum_{u=1}^{r}f^{(u)}(x)(\sum_{v=0}^{u}C_{u}^{v}(-x)^{u-v}\sum_{i=1}^{r}x_{i}^{v}l_{i}(x))\nonumber\\
&&+{\frac
1{(r-1)!}}\sum_{i=1}^{r}l_{i}(x)\int_{x}^{x_{i}}(x_i-s)^{r-1}f^{(r)}(s)ds,
\end{eqnarray*}
which implies that
\begin{eqnarray*}
{\bar{w}(x)}(g(x)-H(g,x))={\frac
1{(r-1)!}}{\bar{w}(x)}\sum_{i=1}^{r}l_{i}(x)\int_{x}^{x_{i}}(x_i-s)^{r-1}g^{(r)}(s)ds,
\end{eqnarray*}
since $|l_{i}(x)|\leqslant C$ for $x\in [x_{r-5/2},x_{r-3/2}],\
i=1,2,\cdots,r$. It follows from
$\frac{|x_i-s|^{r-1}}{{\bar{w}}(s)}\leqslant
\frac{|x_i-x|^{r-1}}{{\bar{w}}(x)},$ $s$ between $x_i$ and $x$, then
\begin{eqnarray*}
{\bar{w}(x)}|g(x)-H(g,x)|&\leqslant&C\bar{w}(x)\sum_{i=1}^{r}\int_{x}^{x_{i}}(x_i-s)^{r-1}|g^{(r)}(s)|ds\nonumber\\
&\leqslant& C\|{\bar{w}}\varphi^{r\lambda}g^{(r)}\|\sum_{i=1}^{r}(x_i-x)^{r-1}\int_{x}^{x_{i}}\varphi^{-r\lambda}(s)ds\nonumber\\
&\leqslant& C{\frac
{\varphi^r(x)}{\varphi^{r\lambda}(x)}}\|{\bar{w}}\varphi^{r\lambda}g^{(r)}\|\sum_{i=1}^{r}(x_i-x)^{r-1}\int_{x}^{x_{i}}\varphi^{-r}(s)ds\\
&\leqslant&C{\frac
{\delta_n^r(x)}{\varphi^{r\lambda}(x)}}\|{\bar{w}}\varphi^{r\lambda}g^{(r)}\|\sum_{i=1}^{r}(x_i-x)^{r-1}\int_{x}^{x_{i}}\varphi^{-r}(s)ds\\
&\leqslant& C(\frac
{\delta_n(x)}{\sqrt{n}\varphi^\lambda(x)})^r\|{\bar{w}}\varphi^{r\lambda}g^{(r)}\|.
\Box
\end{eqnarray*}
\textbf{Lemma 8.} \textit{If $r \in N,\ 0 \leqslant \lambda
\leqslant 1,\ f\in W_{\bar{w},\lambda}^{r},\ \alpha >0,\ $ we have
\begin{eqnarray}
|{\bar{w}(x)}\varphi^{r\lambda}(x){\bar{B}^{(r)}_{n,r-1}(f,x)}|
\leqslant C\|{\bar{w}}\varphi^{r\lambda}f^{(r)}\|.\label{s10}
\end{eqnarray}}
\textbf{Proof.} It follows from $\frac{|t-u|}{{\bar{w}}(u)}\leqslant
\frac{|t-x|}{{\bar{w}}(x)},$ $u$ between $t$ and $x$, let $t=0,$ we
have
\begin{eqnarray*}
n_i^r|\overrightarrow{\Delta}_{\frac 1{n_i}}^{r}{\bar{F}}_{n}{(\frac
k{n_i})}| &=& n_i^r\int_{-\frac {1}{2n_i}}^{\frac {1}{2n_i}} \cdots
\int_{-\frac {1}{2n_i}}^{\frac
{1}{2n_i}}{\bar{F}}_{n}^{(r)}(x+{\frac {rh}{2}}+u_1+ \cdots
+u_r)du_1 \cdots du_r\\
&\leqslant&
Cn_i^r\|{\bar{w}}\varphi^{r\lambda}{\bar{F}_n}^{(r)}\|\int_{-\frac
{1}{2n_i}}^{\frac {1}{2n_i}} \cdots \int_{-\frac {1}{2n_i}}^{\frac
{1}{2n_i}}{\bar{w}}^{-1}(x+{\frac {rh}{2}}+u_1+ \cdots +u_r).\\
&& \varphi^{-r\lambda}(x+{\frac {rh}{2}}+u_1+ \cdots +u_r)du_1
\cdots du_r\\
&=&
Cn_i^r\|{\bar{w}}\varphi^{r\lambda}{\bar{F}_n}^{(r)}\|\int_{-\frac
{1}{2n_i}}^{\frac {1}{2n_i}} \cdots \int_{-\frac {1}{2n_i}}^{\frac
{1}{2n_i}}{\frac {(x+{\frac {rh}{2}}+u_1+ \cdots
+u_r)}{{\bar{w}}(x+{\frac {rh}{2}}+u_1+ \cdots +u_r)}} \cdot \\
&& {\frac {\varphi^{-r\lambda}(x+{\frac {rh}{2}}+u_1+ \cdots
+u_r)}{(x+{\frac {rh}{2}}+u_1+ \cdots +u_r)}}du_1 \cdots du_r\\
&\leqslant&
Cn_i^r\|{\bar{w}}\varphi^{r\lambda}{\bar{F}_n}^{(r)}\|{\frac
{x}{{\bar{w}}(x)}}\int_{-\frac {1}{2n_i}}^{\frac {1}{2n_i}} \cdots
\int_{-\frac {1}{2n_i}}^{\frac {1}{2n_i}}(x+{\frac {rh}{2}}+u_1+
\cdots +u_r)^{-({\frac {r\lambda}{2}}+1)}\cdot\\
&& [1-(x+{\frac {rh}{2}}+u_1+ \cdots +u_r)]^{-{\frac
{r\lambda}{2}}}du_1 \cdots du_r\\
&\leqslant&
C{{\bar{w}}^{-1}(x)}\varphi^{-r\lambda}(x)\|{\bar{w}}\varphi^{r\lambda}{\bar{F}_n}^{(r)}\|\\
&\leqslant&
C{{\bar{w}}^{-1}(x)}\varphi^{-r\lambda}(x)\|{\bar{w}}\varphi^{r\lambda}f^{(r)}\|.
\end{eqnarray*}
By \cite{Totik}, we have
\begin{eqnarray*}
{\bar{B}^{(r)}_{n,r-1}(f,x)}=\sum_{i=0}^{r-2}{\frac
{n_{i}!}{({n_{i}-r})!}}\sum_{k=0}^{n_i-r}C_{i}(n)\overrightarrow{\Delta}_{\frac
1{n_i}}^{r}{\bar{F}}_{n}{(\frac k{n_i})}p_{n_i-r,k}(x).
\end{eqnarray*}
Obviously
\begin{eqnarray*}
|{\bar{w}(x)}\varphi^{r\lambda}(x){\bar{B}^{(r)}_{n,r-1}(f,x)}|
\leqslant C\|{\bar{w}}\varphi^{r\lambda}f^{(r)}\|. \Box
\end{eqnarray*}
\textbf{Lemma 9.} \textit{If $r \in N,\ 0 \leqslant \lambda
\leqslant 1,\ f\in C_{\bar{w}},\ \alpha >0,\ $ we have
\begin{eqnarray}
|{\bar{w}(x)}\varphi^{r\lambda}(x){\bar{B}^{(r)}_{n,r-1}(f,x)}|
\leqslant
Cn^{r/2}\{max\{n^{r(1-\lambda)/2},\varphi^{r(\lambda-1)}\}\}\|{\bar{w}}f\|.\label{s11}
\end{eqnarray}}
\textbf{Proof.} \textit{Case 1.} If $0\leqslant \varphi(x)\leqslant
{\frac {1}{\sqrt{n}}}$, by $(\ref{s4})$, we have
\begin{eqnarray*}
|{\bar{w}(x)}\varphi^{r\lambda}(x){\bar{B}}_{n,r-1}^{(r)}(f,x)|\leqslant
Cn^{-r\lambda/2}|{\bar{w}(x)}{\bar{B}}_{n,r-1}^{(r)}(f,x)|\leqslant
Cn^{r(1-\lambda/2)}\|{\bar{w}}f\|.
\end{eqnarray*} \textit{Case 2.} If
$\varphi(x)> {\frac {1}{\sqrt{n}}}$, we have
\begin{eqnarray*}
&&|{\bar{B}}_{n,r-1}^{(r)}(f,x)|=|B_{n,r-1}^{(r)}({\bar{F}_{n}},x)|\nonumber\\
&\leqslant&(\varphi^{2}(x))^{-r}\sum_{i=0}^{r-2}\sum_{j=0}^{r}Q_{j}(x,n_i)|C_{i}(n)|n_{i}^{j}\sum_{k=0}^{n_i}|(x-{\frac
kn_{i}})^{j}||{\bar{F}}_{n}({\frac kn_{i}})|p_{n_i,k}(x),
\end{eqnarray*}
where\\
$Q_{j}(x,n_i)=(n_ix(1-x))^{[{\frac {r-j}{2}}]},$ and
$(\varphi^{2}(x))^{-r}Q_{j}(x,n_i)n_{i}^{j}\leqslant
C(n_i/\varphi^{2}(x))^{\frac {r+j}{2}}$.\\
So
\begin{eqnarray*}
&&|{\bar{w}(x)}\varphi^{r\lambda}(x){\bar{B}}_{n,r-1}^{(r)}(f,x)|\nonumber\\
&\leqslant&
C{\bar{w}(x)}\varphi^{r\lambda}(x)\sum_{i=0}^{r-2}\sum_{j=0}^{r}({\frac
{n_{i}}{\varphi^2(x)}})^{\frac {r+j}{2}}\sum_{k=0}^{n_i}|(x-{\frac
kn_{i}})^{j}{\bar{F}}_{n}({\frac kn_{i}})|p_{n_i,k}(x)\nonumber\\
&=&
C{\bar{w}(x)}\varphi^{r\lambda}(x)\sum_{i=0}^{r-2}\sum_{j=0}^{r}({\frac
{n_{i}}{\varphi^2(x)}})^{\frac {r+j}{2}}\sum_{k/n_i\in A}|(x-{\frac
kn_{i}})^{j}||{\bar{F}}_{n}({\frac kn_{i}})|p_{n_i,k}(x)\nonumber\\
&&+
C{\bar{w}(x)}\varphi^{r\lambda}(x)\sum_{i=0}^{r-2}\sum_{j=0}^{r}({\frac
{n_{i}}{\varphi^2(x)}})^{\frac {r+j}{2}}\sum_{x_2^{\prime} \leqslant
k/n_i\leqslant x_3^\prime}|{(x-{\frac kn_{i}})^{j}}||H({\frac
kn_{i}})|p_{n_i,k}(x)\nonumber\\
&:=&\sigma_1+ \sigma_2.
\end{eqnarray*}
Where $A:=[0,x_2^{\prime}]\cup [x_3^{\prime},1]$. Working as lemma
3, we can easily get $\sigma_1\leqslant Cn^{\frac
r2}\varphi^{r(\lambda-1)}(x)\|{\bar{w}}f\|,$ and $\sigma_2\leqslant
Cn^{\frac r2}\varphi^{r(\lambda-1)}(x)\|{\bar{w}}f\|.$ By bringing
these facts together, the lemma is proved. $\Box$
%---------------------------------------------------------------------------------------------------------------------------
\section{Proof of Theorem}
\subsection*{The direct theorem}
We know
\begin{eqnarray}
{\bar{F}}_n(t)={\bar{F}}_n(x)+{\bar{F}}'_n(t)(t-x) + \cdots +
{\frac{1}{(r-1)!}}\int_x^t
(t-u)^{r-1}{\bar{F}}^{(r)}_n(u)du,\label{s12}\\
B_{n,r-1}((\cdot-x)^k,x)=0, \ k=1,2,\cdots,r-1.\label{s13}
\end{eqnarray}
According to the definition of $W_{\bar {w},\lambda}^{r},$
\ for any $g \in W_{\bar {w},\lambda}^{r},$ we have
${\bar{B}}_{n,r-1}(g,x)=B_{n,r-1}({\bar{G}}_{n}(g),x),$ and
${\bar{w}(x)}|{\bar{G}}_{n}(x)-B_{n,r-1}({\bar{G}}_{n},x)|={\bar{w}(x)}|B_{n,r-1}(R_r({\bar{G}}_n,t,x),x)|,$
thereof $R_r({\bar{G}}_n,t,x)=\int_x^t
(t-u)^{r-1}{\bar{G}}^{(r)}_n(u)du.$
\\It follows from $\frac{|t-u|^{r-1}}{{\bar{w}}(u)}\leqslant
\frac{|t-x|^{r-1}}{{\bar{w}}(x)},$ $u$ between $t$ and $x$, we have
\begin{eqnarray}
{\bar{w}(x)}|{\bar{G}}_{n}(x)-B_{n,r-1}({\bar{G}}_{n},x)| &
\leqslant &
C\|{\bar{w}}\varphi^{r\lambda}{\bar{G}}^{(r)}_n\|{\bar{w}(x)}B_{n,r-1}(\int_x^t{\frac
{|t-u|^{r-1}}{{\bar{w}(u)}\varphi^{r\lambda}(u)}du,x})\nonumber\\
& \leqslant &
C\|{\bar{w}}\varphi^{r\lambda}{\bar{G}}^{(r)}_n\|{\bar{w}}(x)(B_{n,r-1}(\int_x^t{\frac
{|t-u|^{r-1}}{\varphi^{2r\lambda}(u)}}|du,x))^{\frac 12}\cdot
\nonumber\\
&&(B_{n,r-1}(\int_x^t{\frac
{|t-u|^{r-1}}{{\bar{w}^2(u)}}}du,x))^{\frac 12}.\label{s14}
\end{eqnarray}
also
\begin{eqnarray}
\int_x^t{\frac {|t-u|^{r-1}}{\varphi^{2r\lambda}(u)}}du \leqslant
C{\frac {|t-x|^r}{\varphi^{2r\lambda}(x)}},\ \int_x^t{\frac
{|t-u|^{r-1}}{{\bar{w}^2(u)}}}du \leqslant {\frac
{|t-x|^r}{{\bar{w}^2(x)}}}.\label{s15}
\end{eqnarray}
By (\ref{s3}), (\ref{s14}) and (\ref{s15}), we have
\begin{eqnarray}
{\bar{w}(x)}|{\bar{G}}_{n}(x)-B_{n,r-1}({\bar{G}}_{n},x)|
&\leqslant&
C\|{\bar{w}}\varphi^{r\lambda}{\bar{G}}^{(r)}_n\|\varphi^{-r\lambda}(x)B_{n,r-1}
(|t-x|^r,x)\nonumber\\
&\leqslant& Cn^{-\frac r2}{\frac
{\varphi^{r}(x)}{\varphi^{r\lambda}(x)}}\|{\bar{w}}\varphi^{r\lambda}{\bar{G}}^{(r)}_n\|\nonumber\\
&\leqslant& Cn^{-\frac r2}{\frac
{\delta_n^{r}(x)}{\varphi^{r\lambda}(x)}}\|{\bar{w}}\varphi^{r\lambda}{\bar{G}}^{(r)}_n\|\nonumber\\
&=& C(\frac
{\delta_n(x)}{\sqrt{n}\varphi^\lambda(x)})^r\|{\bar{w}}\varphi^{r\lambda}{\bar{G}}^{(r)}_n\|.\label{s16}
\end{eqnarray}
By (\ref{s9}) and (\ref{s16}), when $g \in
W_{\bar{w},\lambda}^{r},$ then
\begin{eqnarray}
{\bar{w}(x)}|g(x)-{\bar{B}_{n,r-1}(g,x)}| &\leqslant&
{\bar{w}(x)}|g(x)-{\bar{G}}_{n}(g,x)| +
{\bar{w}(x)}|{\bar{G}}_{n}(g,x)-{\bar{B}_{n,r-1}(g,x)}|\nonumber\\
&\leqslant& {\bar{w}(x)}|g(x)-H(g,x)|_{[x_1,x_4]} + C(\frac
{\delta_n(x)}{\sqrt{n}\varphi^\lambda(x)})^r\|{\bar{w}}\varphi^{r\lambda}{\bar{G}}^{(r)}_n\|\nonumber\\
&\leqslant& C(\frac
{\delta_n(x)}{\sqrt{n}\varphi^\lambda(x)})^r\|{\bar{w}}\varphi^{r\lambda}g^{(r)}\|.\label{s17}
\end{eqnarray}
For $f \in C_{\bar{w}},$ we choose proper $g \in W_{\bar
{w},\lambda}^{r},$ by (\ref{s5}) and (\ref{s17}), then
\begin{eqnarray*}
{\bar{w}(x)}|f(x)-{\bar{B}_{n,r-1}(f,x)}| &\leqslant&
{\bar{w}(x)}|f(x)-g(x)| + {\bar{w}(x)}|{\bar{B}_{n,r-1}(f-g,x)}| +
{\bar{w}(x)}|g(x)-{\bar{B}_{n,r-1}(g,x)}|\\
&\leqslant& C(\|{\bar{w}}(f-g)\|+(\frac
{\delta_n(x)}{\sqrt{n}\varphi^\lambda(x)})^r\|{\bar{w}}\varphi^{r\lambda}g^{(r)}\|)\\
&\leqslant& C\omega_{\varphi^\lambda}^r(f,\frac
{\delta_n(x)}{\sqrt{n}\varphi^\lambda(x)})_{\bar{w}}. \Box
\end{eqnarray*}
%------------------------------------------------------------------------------------------------------
\subsection*{The inverse theorem}
We define the weighted main-part modulus for $D=R_+$ by
\begin{eqnarray*}
\Omega_{\varphi^\lambda}^r(C,f,t)_{\bar{w}} = \sup_{0<h \leqslant
t}\|{\bar{w}}\Delta_{{h\varphi}^\lambda}^rf\|_{[Ch^\ast,\infty]},\\
\Omega_{\varphi^\lambda}^r(1,f,t)_{\bar{w}} =
\Omega_{\varphi^\lambda}^r(f,t)_{\bar{w}}.
\end{eqnarray*}
where $C>2^{1/\beta(0)-1},\ \beta(0)>0,$ and $h^\ast$ is given by
\begin{eqnarray*}
h^\ast= \left\{
\begin{array}{lrr}
(Ar)^{1/1-\beta(0)}h^{1/1-\beta(0)},  && 0 \leqslant \beta(0) <1,
    \\
0,   && \beta(0) \geqslant 1.
              \end{array}
\right.
\end{eqnarray*}
The main-part $K$-functional is given by
\begin{eqnarray*}
K_{r,\varphi^\lambda}(f,t^r)_{\bar{w}}=\sup_{0<h \leqslant
t}\inf_g\{\|{\bar{w}}(f-g)\|_{[Ch^\ast,\infty]}+t^r\|{\bar{w}}\varphi^{r\lambda}g^{(r)}\|_{[Ch^\ast,\infty]},\
g^{(r-1)} \in A.C.((Ch^\ast,\infty))\}.
\end{eqnarray*}
By \cite{Totik}, we have
\begin{eqnarray}
C^{-1}\Omega_{\varphi^\lambda}^r(f,t)_{\bar{w}} \leqslant
\omega_{\varphi^\lambda}^{r}(f,t)_{\bar{w}} \leqslant
C\int_0^t{\frac {\Omega_{\varphi^\lambda}^r(f,\tau)_{\bar{w}}}{\tau}}d\tau,\label{s18} \\
C^{-1}K_{r,\varphi^\lambda}(f,t^r)_{\bar{w}} \leqslant
\Omega_{\varphi^\lambda}^r(f,t)_{\bar{w}} \leqslant
CK_{r,\varphi^\lambda}(f,t^r)_{\bar{w}}.\label{s19}
\end{eqnarray}
\textbf{Proof.} Let $\delta>0,$ by (\ref{s19}), we choose proper $g$
so that
\begin{eqnarray}
\|{\bar{w}}(f-g)\| \leqslant
C\Omega_{\varphi^\lambda}^r(f,\delta)_{\bar{w}},\
\|{\bar{w}}\varphi^{r\lambda}g^{(r)}\| \leqslant
C\delta^{-r}\Omega_{\varphi^\lambda}^r(f,\delta)_{\bar{w}}.\label{s20}
\end{eqnarray}
then
\begin{eqnarray*}
|{\bar{w}}(x)\Delta_{h\varphi^\lambda}^rf(x)| &\leqslant&
|{\bar{w}}(x)\Delta_{h\varphi^\lambda}^r(f(x)-{\bar{B}_{n,r-1}(f,x)})|+|{\bar{w}}(x)\Delta_{h\varphi^\lambda}^r\bar{B}_{n,r-1}(f-g,x)|\nonumber\\
&+& |{\bar{w}}(x)\Delta_{h\varphi^\lambda}^r{\bar{B}_{n,r-1}(g,x)}|\nonumber\\
&\leqslant& \sum_{j=0}^rC_r^j(n^{-\frac
12}{\frac {\delta_n(x+({\frac r2}-j)h\varphi^\lambda(x))}{\varphi^\lambda(x+({\frac r2}-j)h\varphi^\lambda(x))}})^{\alpha_0}\nonumber\\
&+& \int_{-{\frac {h\varphi^\lambda(x)}{2}}}^{\frac
{h\varphi^\lambda(x)}{2}}\cdots \int_{-{\frac
{h\varphi^\lambda(x)}{2}}}^{\frac
{h\varphi^\lambda(x)}{2}}{\bar{w}}(x){\bar{B}^{(r)}_{n,r-1}(f-g,x+\sum_{k=1}^ru_k)}du_1\cdots
du_r\nonumber\\
&+& \int_{-{\frac {h\varphi^\lambda(x)}{2}}}^{\frac
{h\varphi^\lambda(x)}{2}}\cdots \int_{-{\frac
{h\varphi^\lambda(x)}{2}}}^{\frac
{h\varphi^\lambda(x)}{2}}{\bar{w}}(x){\bar{B}^{(r)}_{n,r-1}(g,x+\sum_{k=1}^ru_k)}du_1\cdots
du_r\nonumber\\
&:=& J_1+J_2+J_3.
\end{eqnarray*}
Obviously
\begin{eqnarray}
J_1 \leqslant C(n^{-\frac 12}\varphi^{-\lambda}(x)\delta_n(x))^{\alpha_0}.\label{s21}
\end{eqnarray}
By (\ref{s4}) and (\ref{s20}), we have
\begin{eqnarray}
J_2 &\leqslant& Cn^r\|{\bar{w}}(f-g)\|\int_{-{\frac
{h\varphi^\lambda(x)}{2}}}^{\frac {h\varphi^\lambda(x)}{2}}\cdots
\int_{-{\frac {h\varphi^\lambda(x)}{2}}}^{\frac
{h\varphi^\lambda(x)}{2}}du_1 \cdots du_r\nonumber\\
&\leqslant& Cn^rh^r\varphi^{r\lambda}(x)\|{\bar{w}}(f-g)\|\nonumber\\
&\leqslant&
Cn^rh^r\varphi^{r\lambda}(x)\Omega_{\varphi^\lambda}^r(f,\delta)_{\bar{w}}.\label{s22}
\end{eqnarray}
By (\ref{s11}), we let $\lambda=1,$ and (\ref{s6}) as well
as (\ref{s20}), we have
\begin{eqnarray}
J_2 &\leqslant& Cn^{\frac r2}\|{\bar{w}}(f-g)\|\int_{-{\frac
{h\varphi^\lambda(x)}{2}}}^{\frac {h\varphi^\lambda(x)}{2}} \cdots
\int_{-{\frac {h\varphi^\lambda(x)}{2}}}^{\frac
{h\varphi^\lambda(x)}{2}}\varphi^{-r}(x+\sum_{k=1}^ru_k)du_1 \cdots du_r\nonumber\\
&\leqslant& Cn^{\frac r2}h^r\varphi^{r(\lambda-1)}(x)\|{\bar{w}}(f-g)\|\nonumber\\
&\leqslant& Cn^{\frac
r2}h^r\varphi^{r(\lambda-1)}(x)\Omega_{\varphi^\lambda}^r(f,\delta)_{\bar{w}}.\label{s23}
\end{eqnarray}
By (\ref{s10}) and (\ref{s20}), we have
\begin{eqnarray}
J_3 &\leqslant&
C\|{\bar{w}}\varphi^{r\lambda}g^{(r)}\|{\bar{w}(x)}\int_{-{\frac
{h\varphi^\lambda(x)}{2}}}^{\frac {h\varphi^\lambda(x)}{2}} \cdots
\int_{-{\frac {h\varphi^\lambda(x)}{2}}}^{\frac
{h\varphi^\lambda(x)}{2}}{\bar{w}^{-1}(x+\sum_{k=1}^ru_k)}\varphi^{-r\lambda}(x+\sum_{k=1}^ru_k)du_1 \cdots du_r\nonumber\\
&\leqslant& Ch^r\|{\bar{w}}\varphi^{r\lambda}g^{(r)}\|\nonumber\\
&\leqslant&
Ch^r\delta^{-r}\Omega_{\varphi^\lambda}^r(f,\delta)_{\bar{w}}.\label{s24}
\end{eqnarray}
Now, by (\ref{s21}), (\ref{s22}), (\ref{s23}) and
(\ref{s24}), we get
\begin{eqnarray*}
|{\bar{w}}(x)\Delta_{h\varphi^\lambda}^rf(x)| &\leqslant&
C\{(n^{-\frac 12}\delta_n(x))^{\alpha_0} + h^r(n^{-\frac
12}\delta_n(x))^{-r}\Omega_{\varphi^\lambda}^r(f,\delta)_{\bar{w}} +
h^r\delta^{-r}\Omega_{\varphi^\lambda}^r(f,\delta)_{\bar{w}}\}.
\end{eqnarray*}
When $n \geqslant 2,$ we have
\begin{eqnarray*}
n^{-\frac 12}\delta_n(x) < (n-1)^{-\frac 12}\delta_{n-1}(x)
\leqslant \sqrt{2}n^{-\frac 12}\delta_n(x),
\end{eqnarray*}
Choosing proper $x, n \in N,$ so that
\begin{eqnarray*}
n^{-\frac 12}\delta_n(x) \leqslant \delta < (n-1)^{-\frac
12}\delta_{n-1}(x),
\end{eqnarray*}
Therefore
\begin{eqnarray*}
|{\bar{w}}(x)\Delta_{h\varphi^\lambda}^rf(x)| \leqslant
C\{\delta^{\alpha_0} +
h^r\delta^{-r}\Omega_{\varphi^\lambda}^r(f,\delta)_{\bar{w}}\}.
\end{eqnarray*}
By Borens-Lorentz lemma, we get
\begin{eqnarray}
\Omega_{\varphi^\lambda}^r(f,t)_{\bar{w}} \leqslant
Ct^{\alpha_0}.\label{s25}
\end{eqnarray}
So, by (\ref{s25}), we get
\begin{eqnarray*}
\omega_{\varphi^\lambda}^{r}(f,t)_{\bar{w}} \leqslant
C\int_0^t{\frac
{\Omega_{\varphi^\lambda}^r(f,\tau)_{\bar{w}}}{\tau}}d\tau =
C\int_0^t\tau^{\alpha_0-1}d\tau = Ct^{\alpha_0}. \Box
\end{eqnarray*}
%---------------------------------------------------------------------------------------------------------------------------------------------

\end{document}